\documentclass[12pt, a4paper]{article}
\usepackage{fullpage}
\usepackage[top=2cm, bottom=2.5cm, left=2.5cm, right=2.5cm]{geometry}
\usepackage{amsmath,amsthm,amsfonts,amssymb,amscd}
\usepackage{comment}
\usepackage{lastpage}
\usepackage{enumerate}
\usepackage{fancyhdr}
\usepackage{mathrsfs}
\usepackage{xcolor}
\usepackage{graphicx}
\usepackage{listings}
\usepackage{hyperref}
\usepackage{mdframed}
\usepackage{amssymb}
\usepackage{amsthm}
\usepackage{mathtools}
\usepackage[square,sort,comma,numbers]{natbib}
\usepackage{wrapfig}
\usepackage{textcomp}
\usepackage{graphicx}
\usepackage{float}
\usepackage[labelformat=empty]{caption}
\usepackage{hyperref}% http://ctan.org/pkg/hyperref
\usepackage{cleveref}% http://ctan.org/pkg/cleveref
\usepackage[makeroom]{cancel}
\usepackage{xpatch}
\usepackage{parskip}
\allowdisplaybreaks

\usepackage{listofitems}
\newcommand\inner[1]{\readlist\innerlist{#1}\langle%
  \ifnum\innerlistlen=1\relax#1,#1\else#1\fi\rangle}

\newcommand\covariance[1]{\readlist\innerlist{#1}\text{Cov}\left(%
  \ifnum\innerlistlen=1\relax#1,#1\else#1\fi\right)}

\usepackage{mathtools}

\newcommand\nnfootnote[1]{%
  \begin{NoHyper}
  \renewcommand\thefootnote{}\footnote{#1}%
  \addtocounter{footnote}{-1}%
  \end{NoHyper}
}

\def\valitemsep{-0.3em} % itemize
\title{Mathematical Artifacts Have Politics: \\ The Journey from Examples to Embedded Ethics}
\author{
Dennis M\"uller\footnote{RWTH Aachen University, Germany. \texttt{\href{mailto:dennis.mueller3@rwth-aachen.de}{dennis.mueller3@rwth-aachen.de}}.}  \\ Maurice Chiodo\footnote{Centre for the Study of Existential Risk, University of Cambridge, United Kingdom. \texttt{\href{mailto:mcc56@cam.ac.uk}{mcc56@cam.ac.uk}}.}    
}

\begin{document}
\maketitle
\begin{abstract}
We extend Langdon Winner's idea that artifacts have politics into the realm of mathematics. To do so, we first provide a list of examples showing the \textit{existence} of mathematical artifacts that have politics. In the second step, we provide an argument that shows that \textit{all} mathematical artifacts have politics. We conclude by showing the implications for embedding ethics into mathematical curricula. We show how acknowledging that mathematical artifacts have politics can help mathematicians design better exercises for their mathematics students. 
\end{abstract}
\noindent \textit{Keywords:} artifacts have politics, embedded ethics, mathematics in society, ethics in mathematics, mathematics education
\nnfootnote{2020 \textit{AMS Classification:} 01A80. 00A30, 97A40, 97D20.}
\setcounter{tocdepth}{1}
\tableofcontents

\newpage
\section{Introduction}
In this essay, we present a three-step argument that both supports the idea of embedded ethics in mathematics, and refutes the concept of political neutrality of mathematical knowledge. Our first step, providing a list of examples demonstrating how a modern society can use and be affected by mathematics, aims to prove the \textit{existence} of mathematical artifacts that have politics. These examples predicate our second step, in which we provide an argument for why \textit{every} mathematical artifact has politics. This then provides the \textit{validation} and insights for the final step of embedding ethics into mathematical teaching and practice. In doing so, we extend Winner's famous essay \textit{Do artifacts have politics?} focusing on engineering artifacts\cite{Winner1980}, and the many works analysing the politics of artificial intelligence 
(e.g. \cite{inglesdoes, lazovich2020does, winner1978autonomous}), into the world of mathematics.

\section{Some Uses and Effects of Mathematics in Society}
\label{sec:useandeffects}
In this section we list and describe \textit{some} of the uses and effects of mathematics in society. Of course, no such list can ever be complete because new applications of mathematics are found daily. Nonetheless, it is necessary to see a selection of uses and effects in some detail to get a better feeling for what it means for an artifact to have politics, before we explore a more general argument in section \ref{sec:all_artifacts_have_politics} and to better understand the elements of embedded ethics in section \ref{sec:embedded_ethics}. 

Such a list, even if incomplete, can serve as a useful starting point and tool for creating interesting and challenging mathematical exercises for students because it identifies some of the underlying dynamics involved in the general use of mathematics. Any area of mathematics can be connected to some use or effect, and as such, these exemplary uses can help us to design mathematical exercises that explore mathematics and its politics simultaneously.

\begin{figure}[ht!]
    \centering
    \includegraphics[width=0.9\textwidth]{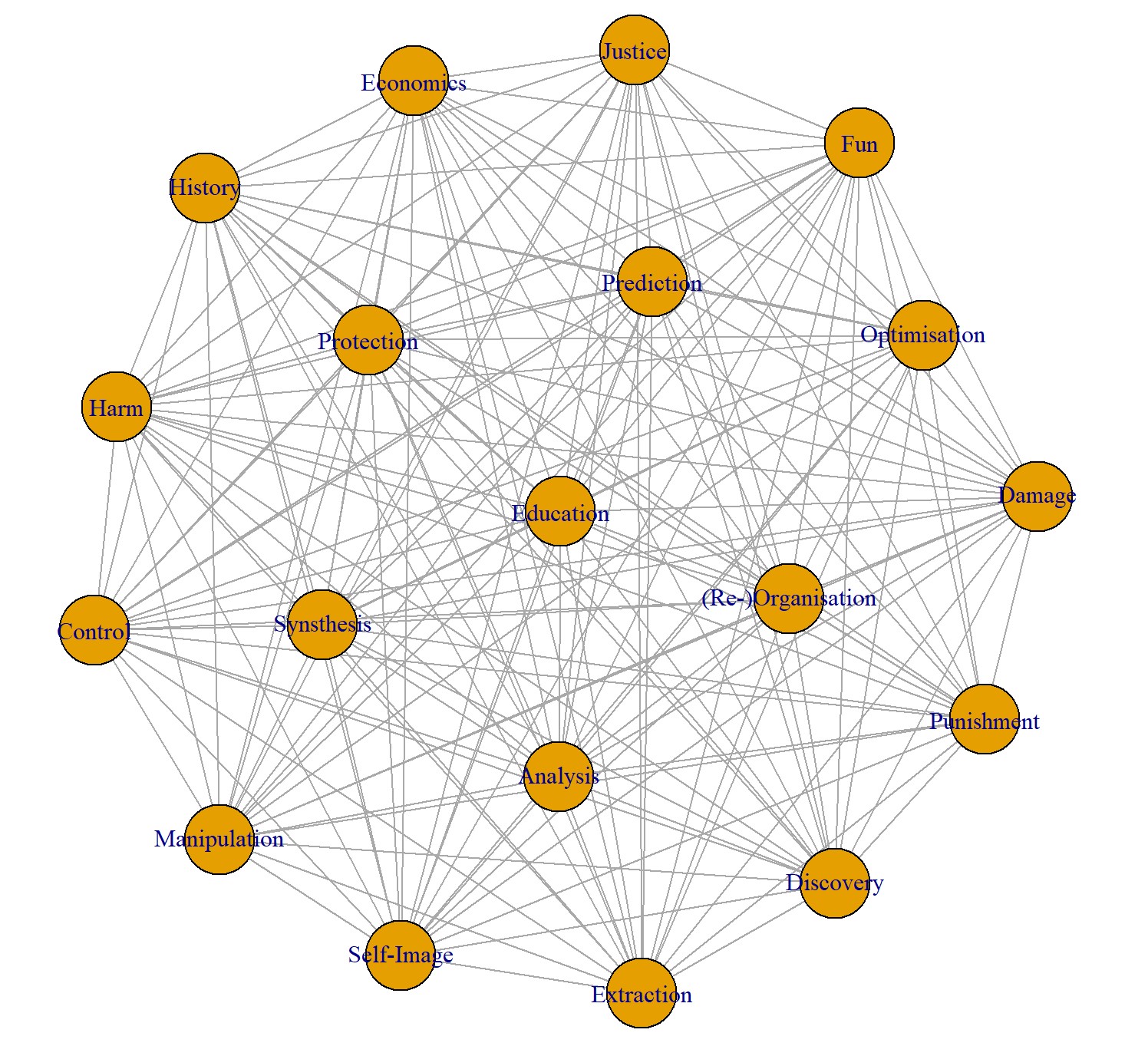}
    \caption{Figure 1: All uses and effects are connected to each other.}
    \label{fig:fully_connected_graph}
\end{figure}

The uses and effects we will cover are: prediction, discovery, (re-)organisation, optimisation, protection, fun, extraction, analysis, synthesis, damage, harm, education, punishment, justice, control, self-image, economics, history, and manipulation; see Figure \ref{fig:fully_connected_graph}. 

\subsection{Prediction}
Predictions typically include statistical results to infer from a sample to an overall population, or a quantitative scientific theory about nature or society to predict the future or an otherwise unknown socio-technical, technical, biological or physical feature. Any prediction alters the state of the world, even when it is not acted on, as it increases our knowledge and potential for actions. When acted on, predictions 
\begin{itemize} 
\setlength\itemsep{\valitemsep}
    \item about systems involving human interaction can become a self-fulfilling prophecy and self-destructive as they can amplify or reduce certain behaviours (cf. \cite{Boudry08062023}),
    \item help us to imagine and navigate an unknown future, and by contrast, assist us in understanding the past (cf. \cite{beckert2016imagined}),
    \item can change the minds of those who perform or use them on a biological and psychological level (cf. \cite{HutchinsonBarrettBrain}). 
\end{itemize}
\noindent\textbf{Examples:} Feedback loops in recommender machine learning systems can lead to echo chambers\cite{jiang2019degenerate}; Covid-19 models have driven policy decisions and put mathematical models and their makers into the spotlight of political attention\cite{mcbryde2020role, questionsofresponsibility}, and so have models about climate change\cite{lemos2010climate}; search engines try to predict the most relevant result for their users, but when they are systematically wrong, they can induce forms of oppression, leading to lasting effects on minority groups\cite{noble2018algorithms}; Einstein's general relativity predicts the existence of gravitational waves\cite{cervantes2016brief}; chaos theory gives ways to understand and potentially overcome some of the limits of reliability and predictability\cite{boccaletti2000control}; mathematical predictions play an essential role in creating imaginations and expectations of the future, and thus drive much of our individual and collective economic activity\cite{beckert2016imagined}.

\subsection{Discovery}
As part of the sciences, politics or everyday life, we use mathematics in an attempt to discover truths about the physical and socio-technical world surrounding us. In the mode of discovery, we fundamentally build on 
\begin{itemize}
\setlength\itemsep{\valitemsep}
    \item the trustworthiness and deductive power of mathematical results to establish trustworthy theories of the world (cf. \cite{Morton1933, popper2005logic}),
    \item the power of mathematics to construct search spaces (cf. \cite{spirtes2000causation}), 
    \item and abstraction to simplify the object of discovery for it to be captured and formalised by mathematical knowledge (cf. \cite{rosenblueth_wiener_1945, Thompson.2022, frantz1995taxonomy}).
\end{itemize}
\noindent\textbf{Examples:} Predictive policing and crime forecasts can support the everyday job of police services\cite{perry2013predictive}, but they also lead to deep political questions about fairness, bias and privacy\cite{alikhademi2022review}; mathematics underlies the search engines that organise the world's information and thus has become fundamental to the discovery of knowledge in everyone's lives\cite{brin1998anatomy}, thereby deeply affecting how people see the world and understand themselves\cite{noble2018algorithms}; mathematical physics helped to search for, and build the machines detecting gravitational waves\cite{cervantes2016brief}; models for climate change cannot capture everything, and the abstraction required often leads to deep philosophical and methodological debates about their construction and use\cite{frigg2015philosophy}.

\subsection{(Re-)Organisation}
Mathematics can be used to (re-)organise technical or socio-technical systems. Such usage 
\begin{itemize}
\setlength\itemsep{\valitemsep}
    \item often builds on potentially conflicting metrics, indicators and measures on which the (re-)organisation is judged (cf. \cite{muller2018tyranny}), 
    \item regularly requires trade-offs and choices, for example, between importance and sensitivity (cf. \cite{da2001exploring}),
    \item and often is within the framework of, or breaks with, a widely accepted scientific theory\cite{kuhn2012structure}, or may break with their practitioner's expectations (cf. \cite{franccois2010revolutions, gillies1992revolutions}) or other's expectations.
\end{itemize}
\noindent\textbf{Examples:} Different fairness measures can be contradictory\cite{hutchinson201950, friedler2021possibility}; machine learning algorithms may require a trade-off between fairness and interpretability to maintain accurate predictions\cite{agarwal2021trade}; many of the decision problems used in operations research (e.g. mixed integer programming) are NP-complete, while others (such as linear programming\cite{hofman2006linear}) can be solved in polynomial time, leading to trade-offs in model selection; there exists an ever-growing list of counterintuitive results in basic probability\cite{nickerson2007penney}; the experimental observations along with new mathematical models broke with the Aristotelian worldview and led to the scientific revolution in the early modern world\cite{shapin2018scientific}.

\subsection{Optimisation}
Any form of optimisation builds on deep ideas about how the world is organised, and how it should progress and develop. Some of those are forms of colonial or scientific knowledge that have become standard in the industrialised modern world\citep[p.~97]{mckelveyneves}. Through this, optimisation can manifest itself in many ways, including 
\begin{itemize}
\setlength\itemsep{\valitemsep}
    \item as mathematical artifacts building on quantification, abstraction, formalism, generalisation, idealisation and other forms of de-situationing from the specificities of the problem (cf. \cite{mckelveyneves}),
    \item as a human practice building on people's aspirations, emotions, desires, imaginations and needs (cf. \cite{mckelveyneves, beckert2016imagined, Ziewitz.2019}),
    \item and as a tool of power and legitimisation it can provide us with answers to the question of how a physical, technical or socio-technical system, organisation or order can, or should, look like (cf. \cite{mckelveyneves, basdevant2023variational}).
\end{itemize}
\noindent\textbf{Examples:} Proxy variables in the measurement of our economies, such as those used to model living standards\cite{Montgomery.2000}, can become distanced from people's lives and thus may struggle to capture reality; under the hood of modern search engines lie difficult optimisation problems, and as systems, they build on their creator's aspirations, imaginations and technical needs (cf. the initial construction of the Google search engine\cite{brin1998anatomy}); as tools of power, the optimisation techniques employed within a management or operations research context can overstep and change social barriers in various ways\cite{roy1981optimisation}; the variational principles underlying much of physics are an example of how optimisation appears in the physical world\cite{basdevant2023variational}. 

Optimisation often crosses the lines between the factual and normative, e.g. the statistical thinking and optimisation ideas behind Darwin's work on natural evolution\cite{darwin_2009, ruse2009cambridge} have led to forms of Social Darwinism\cite{hawkins1997social}, which then manifested itself in normative ideas about an optimal world.

\subsection{Protection} 
We can use mathematics to protect socio-technical or biological systems from each other. In doing so, we 
\begin{itemize}
\setlength\itemsep{\valitemsep}
    \item often aim to create a mathematical system that is aligned with our norms, values, utilities and laws (cf. \cite{yudkowsky2016ai, hou2023multi, asimov2004robot}),
    \item may build on the epistemological supremacy, techniques and methodology of pure mathematics to build (hopefully) provably secure layers of technical protection for people, institutions, processes or data (cf. \cite{koblitz2007another, rogaway2015moral}),
    \item use mathematical methods to understand the evolution of biological systems and their struggle for existence, in order to advise on their preservation or destruction (cf. \cite{gause2019struggle}).
\end{itemize}
\noindent\textbf{Examples:} AI alignment research tries to understand how to put ethical constraints into socio-technical systems\cite{jiang2019degenerate} and, more generally, how to align AI systems with our human values\cite{yudkowsky2016ai}; many hope that predictive policing\cite{perry2013predictive} will make our police services more efficient, and thus, save lives; selecting interpretable machine learning techniques\cite{molnar2020interpretable} over uninterpretable, potentially more accurate models, can be a means to protect society from the systems it builds; as a means of protection, cryptography has become essential in all our lives and modern infrastructure\cite{rogaway2015moral}; mathematical surveillance systems built to protect us can lead to challenging questions when their existence and scale comes to light\cite{verble2014nsa}.

\subsection{Fun}
We can use mathematics to derive enjoyment. This may involve 
\begin{itemize}
\setlength\itemsep{\valitemsep}
    \item seeing mathematics practice as a puzzle, or solving and constructing mathematical puzzles (cf. \cite{winkler2003mathematical,nova}),
    \item appreciating its beauty and understanding it as a form of art (cf. \cite{mcallister2005mathematical,cellucci2015mathematical,zeki2014experience}),
    \item the social activity of doing mathematics in a group of like-minded people or by simply being good at it (cf. \cite{papanastasiou2004math}).
\end{itemize}
\noindent\textbf{Examples:} Many popular books on mathematical puzzles (e.g. \cite{winkler2003mathematical}), beautiful proofs (e.g. \cite{aigner2015buch}) and beautiful mathematical identities (e.g. \cite{erickson2011beautiful}) exist on the market; the beauty and puzzling nature of mathematics can become part of one's definition of doing mathematics, thereby setting the foundations for discussions around ethics in mathematics\cite{MüllerDP2}.

\subsection{Extraction}
Mathematics can be used to empower the extraction of natural, social or economic resources, or its use may require the extraction of said resources. In doing so, it
\begin{itemize}
\setlength\itemsep{\valitemsep}
    \item may have a morally, legally, ecologically, socially or economically (un-)acceptable footprint (cf. \cite{chiodo2023manifesto, crawford2021atlas}),
    \item can directly or indirectly reduce or increase the consumption of resources by other institutions, groups or systems (cf. \cite{chiodo2023manifesto, crawford2021atlas})
    \item or lead to other planetary and social costs and imbalances. The extraction of natural resources may lead to advantages for some to the detriment of others (cf. \cite{chiodo2023manifesto, crawford2021atlas}).
\end{itemize}
\noindent\textbf{Examples:} The planetary and social footprint of modern AI systems is growing massively and is no longer negligible\cite{crawford2021atlas}; Bitcoin has a famously large energy consumption, and in order to be profitable, mining often requires specialist hardware\cite{o2014bitcoin, vranken2017sustainability}.

\subsection{Analysis and Synthesis}
Mathematics gives us the logical and quantitative tools to 
\begin{itemize}
\setlength\itemsep{\valitemsep}
    \item break down mathematics itself and much of the world surrounding us (cf. \cite{otte1997analysis, Wigner1960Unresonable}),
    \item create new abstract, physical or digital items, systems or structures (cf. \cite{otte1997analysis, Wigner1960Unresonable}),
    \item and to deploy new creations efficiently by amplifying their reach (cf. \cite{muller2022hippocratic, chiodo2023manifesto}).
\end{itemize}
\noindent\textbf{Examples: } The effectiveness and amplification power of mathematics in Big Data has led to a new form of capitalism, nowadays often called ``Surveillance Capitalism''\cite{zuboff2023age}; modern science would not be possible without modern mathematics and the effectiveness of mathematical forms of reasoning, quantification and argumentation (cf. \cite{Wigner1960Unresonable}). 

\subsection{Damage or Harm}
Mathematics can be accidentally, willfully or unknowingly used to do damage or harm. Such actions may include
\begin{itemize}
\setlength\itemsep{\valitemsep}
    \item abusing the power and standing of mathematical arguments (cf. \cite{chiodo2022mathematicians}),
    \item the act of learning mathematics (cf. \cite{ernest2018ethics}),
    \item and preventing it may require a holistic look at the target of the mathematics, its social context and a mathematician's training (cf. \cite{muller2022hippocratic}).
\end{itemize}
\noindent\textbf{Examples:} Due to being hard to challenge, statistical reasoning has misled judges to order false convictions\cite{nobles2005misleading}, potentially requiring professional bodies to stand up and take a position in public (the Royal Statistical Society issued a public statement in the case of R v Sally Clark \cite{RSS-sally-clarke}); the misuses and abuses of mathematics are almost never-ending, and thus can be called ``weapons of maths destruction'' endangering our modern democratic societies\cite{o2017weapons}; AI systems can learn to be racist or antisemitic from biased data or harmful user interactions\cite{handelman2022artificial}; the potential for misuse and harm has led to many calls for Hippocratic oaths (for an overview, see \cite{muller2022hippocratic}).

\subsection{Education}
Mathematics education is deeply connected to the organisation and structure of modern societies. At all levels, mathematics education
\begin{itemize}
\setlength\itemsep{\valitemsep}
    \item is not blind to existing societal problems, such as issues of racial discrimination, equity and fairness (cf. \cite{botelho2015racial, lubienski2002research, civil2007building, forgasz2012towards,martin2019equity,martin2010mathematics}),
    \item can promote or reduce a student's self-worth, anxiety and feelings of accomplishment (cf. \cite{pajares1995mathematics, lawrence2009link, cano2018students}),
    \item can be hindered or fostered by one's attitude towards mathematics (cf. \cite{nicolaidou2003attitudes}).
\end{itemize}
\noindent\textbf{Examples:} The discourse surrounding mathematics education may include discussions about equity, social responsibility and specialised ethics frameworks (e.g. \cite{rycroftsmith2022useful}); some have called good mathematics education a civil right because of its necessity for people to function in and successfully navigate modern society\cite{moses2002radical, buell2019intro}. 

\subsection{Punishment and Justice}
Mathematics can be used to falsely or correctly support, execute or promote punishment
\begin{itemize}
\setlength\itemsep{\valitemsep}
    \item in the legal system through correct or incorrect mathematical arguments (cf. \cite{leung2002prosecutor}),
    \item in education through the use of problematic educational philosophies (cf. \cite{ernest2018ethics, rycroftsmith2022useful}),
    \item or in other institutional settings through the implementation of problematic performance metrics or other quantitative indicators (cf. \cite{muller2018tyranny}).
\end{itemize}
\noindent\textbf{Examples:} Misleading statistics have led to unjustified convictions\cite{nobles2005misleading}; the use of predictive algorithms for A-level grades created a public outcry due to its potential for injustices\cite{heaton2023algorithm}; the overuse of metric-based decision-making can be punishing to those who are not adequately captured by them\cite{muller2018tyranny}.

\subsection{Control}
We can use mathematics to exert control over biological, technical, digital or physical systems, social groups or institutions. In doing so, we 
\begin{itemize}
\setlength\itemsep{\valitemsep}
\item have created the matured mathematical area of control theory that developed many of its modern foundations in the political circumstances of war (cf. \cite{PeschPlail2012}), 
\item may use mathematics to help organise, classify or protect certain aspects or groups of our modern societies (cf. \cite{ferguson2016policing,meijer2019predictive, gates2020school}),
\item have used, misused and abused it in various situations (cf. \cite{wilmott2000use}).
\end{itemize}
\noindent\textbf{Examples:} The transport problem can be used to create and destroy rail networks\cite{Schrijver.2002}; the form, act and results of grading students may be impacted by a student's socio-economic background and social capital(cf. \cite{hjorth2022grading}); the Gaussian copula, used in financial models that estimated and attempted to control risks, played a fundamental role in the global financial crisis of 2008\cite{mackenzie2014formula, salmon2009recipe}.

\subsection{Self-Image}
The existence and use of modern mathematics, and its deep roots in the Greco-Roman tradition of thought, impact the self-image of
\begin{itemize}
\setlength\itemsep{\valitemsep}
    \item all people through its dominance in modern education and everyday life (cf. \cite{ernest2018ethics}), 
    \item all our institutions through the promotion of quantifiable decision-making and other rational forms of reasoning (cf. \cite{tasic2001mathematics, walshaw2004mathematics, muller2018tyranny},
    \item and other societies worldwide throughout history and in the present (cf. \cite{D’Ambrosio2016, d1985ethnomathematics}).
\end{itemize}
\noindent\textbf{Examples:} There is a heated debate about whether one must, should or can decolonise mathematics (e.g.\cite{garcia2022decolonizing, raju2018decolonising, borovik2023decolonisation}). These debates, as well as similarly situated discourses on ethics in mathematics and mathematics for social justice, regularly challenge or defend the potential universality of mathematics and other self-understandings of the discipline\citep[p.~8]{müller2022situating}. 

\subsection{Economics}
Building on its solid foundations and modes of reasoning, mathematics has impacted modern capitalist societies at various levels, e.g. mathematics 
\begin{itemize}
\setlength\itemsep{\valitemsep}
    \item has revolutionised economics helping it become a quantitative rather than purely social science (cf. \cite{screpanti2005outline})
    \item is the cornerstone of many of the organisational principles of modern economies and capitalist dynamics (cf. \cite{beckert2016imagined}),
    \item can be named a production factor of the modern digital economies, right next to land, labour, capital and entrepreneur (cf. \cite{grotschel2010production}).
\end{itemize}
\noindent\textbf{Examples:} The introduction of utility functions has led to a ``marginal revolution'' in economics\cite{moscati2018measuring}, clashing with approaches focusing on the social psychology of human behaviour (cf. \cite{akerlof2010animal, fox2009myth}), which can see economics as a predominantly social and not mathematical science which is (thus) unable to always adequately quantify human decision making; standard measures of economics (e.g. the GPD\cite{england1998measurement}) are used to drive policy, even though they regularly fail to capture social or environmental well-being; the efficiency and amplification potential of mathematical digital products and services have lead to a form of ``surveillance capitalism''\cite{zuboff2023age}.

\subsection{History}
The history of mathematics is in some aspects very similar to the history of other natural sciences, i.e.
\begin{itemize}
\setlength\itemsep{\valitemsep}
    \item it is a history of many different forms and cultures engaging in mathematics (cf. \cite{D’Ambrosio2016, d1985ethnomathematics, boyer2011history}),
    \item it cannot be written as a history of linear, constant progress (cf. \cite{edgerton2004linear, gray2004anxiety}),
    \item is full of different reasons for doing mathematics (cf. \cite{bursill2002we}).
\end{itemize}
\noindent\textbf{Examples:} The study of ethnomathematics has led to new challenges in the writing of histories of mathematics (cf. \cite{D’Ambrosio2016}); geometry was done for different reasons, with different methods, and different intentions throughout the ages\cite{bursill2002we}; the 19th century may have been a century of anxiety for many mathematicians\cite{gray2004anxiety}.

\subsection{Manipulation}
Rational quantification through mathematics regularly sets the standard for proper reasoning about the world (cf. \cite{popper2005logic, descartes2012discourse}). In this context, mathematics can be used 
\begin{itemize}
\setlength\itemsep{\valitemsep}
    \item to efficiently and effectively manipulate the physical world (cf. \cite{Wigner1960Unresonable, Thompson.2022}),
    \item to (in)effectively manipulate the social world using potentially unfalsifiable quantitative theories and models (cf. \cite{Velupillai2005ineffectiveness, mckelveyneves, Thompson.2022}),
    \item and to manipulate the world of socio-technical systems by embedding one's politics into them (cf. \cite{Morris.2021, Thompson.2022}).
\end{itemize}
\noindent\textbf{Examples:} Modern physics and engineering would not be possible without modern mathematics, and it is ``unreasonably effective''\cite{Wigner1960Unresonable}; the machine learning models behind social media's ranking algorithms can promote certain beauty standards and thus inadvertently manipulate viewers and content producers(e.g. \cite{bishop2018anxiety}); micro-targeting as part of political advertising, potentially crossing the line between persuasion and manipulation\citep[p.~86]{berghel2018malice}.

\section{Mathematical Artifacts Have Politics} \label{sec:all_artifacts_have_politics}
In the previous section, we provided examples of how mathematical artifacts can have a political dimension, thereby establishing the existence of \textit{some} mathematical artifacts that have politics. This section will argue that \textit{all} mathematical artifacts have politics, just like Langdon Winner\cite{Winner1980} argued that all technical artifacts have politics. This has also been asserted in the recent \textit{Manifesto for the Responsible Development for Mathematical Works}, albeit without a complete argument\cite{chiodo2023manifesto}. Of course, any such argument must lie on certain assumptions. Ours builds on two. We take as given the lessons from critical theory (cf. \cite{bohman2005critical}), i.e.
\begin{itemize}
    \item \textit{Politics is any action involving or affecting you or someone else directly or indirectly, i.e. them as a person or their values, norms, wishes or desires.} 
    \item \textit{Every activity which is exercised by a human or by a machine is political. Any activity always affects or is affected by you or others, as it is supported or opposed by you or others or it is related to your or someone else's values, norms, wishes or desires.}
\end{itemize}  
We understand this as \textit{morally necessary}, since otherwise there exist moral vacuums in which people can commit harm without the moral judgement of themselves or others. Implicit in this moral necessity is a form of Kant's categorical imperative, i.e. every individual must not solely be treated as a means to an end but always be treated as an end in itself\cite{kant_2011}. 

Why is this relevant for mathematics? A unique freedom governs mathematical practice. The structure of a mathematical argument might look inevitable, and every line logically follows from its predecessor, yet it is full of human decisions and choices\cite{nickel2005ethik}. These decisions reflect the mathematician's or someone else's\footnote{These can include your students, colleagues, your funder, manager, your users, and many other affected parties, see also \cite{chiodo2023manifesto}.} desires, needs, wishes or norms. As such, every mathematical decision or choice is political and contains deep normative and political questions (this is especially true about foundational mathematical decisions, see \cite{wagner2023ethical}). Consequently, every mathematical artifact, i.e. the product of a human or machine performing a mathematical task, has politics. 

Yes, this means that ``2+3=5'' is a mathematical artifact that has politics. Without further context, it does not imply that it is good or bad, but it is not politically neutral. The politics of such a statement include, but are not limited to, the use of widely accepted notations, the use of a widely accepted numbering and counting system, as well as the normative component of writing down a true equation within that system instead of writing it down in alternative systems found throughout history or in indigenous communities around the world\cite{owens2020revising}. Even simple acts like doing arithmetic or counting are social and political\cite{bier2017quantification}. Just like the (natural) sciences' success and progress are intimately connected to its norms and standards\cite[p.~8]{sismondo2010introduction}, so is that of mathematics. The meaning and context of what even simple equations like ``2+3=5'' stand for are not politically neutral. 

Thus, all mathematical artifacts are material actors in our world (cf. \cite{latour1990actor}) that have politics beyond the politics of their maker. They can promote or restrict certain actions, styles of thought or modes of reasoning, and what feels natural to us (e.g., the construction of natural numbers using Peano's axioms\cite{segre1994peano}), might not be all too natural for someone who did not grow up with our heavily axiomatised form of arithmetic using Hindu-Arabic numerals.

\begin{mdframed}
Consider, as another example, the Neutron Diffusion Equation governing the diffusion of neutrons in different materials:
    \[
    \frac{\partial \Phi}{\partial t}= \frac{\sigma-1}{\tau}\Phi + \frac{\lambda^{2} }{3 \tau}\nabla^{2} \Phi
    \]
where $\Phi(\mathbf{x},t)$ is the (free) neutron density, $\sigma$ is the average number of neutrons released in a fission event, $\tau$ is the average time between fission events, and $\lambda$ is the average distance a neutron travels before being absorbed by a nucleus. This equation, and its mathematical solution showing that the mass of Uranium needed for an atomic bomb is a feasibly obtainable quantity, is an example where knowledge of the existence had huge political influence before anything remotely similar to a bomb was even built with it. In a very real sense, it was enough to shift the war effort during World War 2\cite{serber1992alamos}. The equation had ambition deeply embedded into it, giving us both the atomic bomb and nuclear power plants, but its politics go beyond dual-use issues of research with its varying forms of political and ethical interpretations and points of view.   
\end{mdframed}

The definition of politics also means that the decision to create a mathematical output without obvious applications is deeply political as it follows someone's values, norms, wishes or desires; depending on its motivation and the consciousness of the decision. In particular, decisions not to be political (with or about one's or another's mathematics) are always political. Once again, this does not imply that it is good or bad. It can be either, depending on the specific circumstances of these decisions. But politically neutral, it is not. 

Hence, in some sense, the statement that ``mathematical artifacts have politics'' is, at the same time, extremely benign and radical. 

Unlike much of the modern political discourse that we are used to, the statement that ``mathematical artifacts have politics'' neither claims the moral ground of good or bad, nor tells 
you how to believe or act. The statement is thus fundamentally different from many of the notions of politics that are found in the news, social media or other forms of discourse. In essence, ``mathematical artifacts have politics'' does not mean they have ``one kind of politics''. Each artifact can come with a spectrum of politics that depends on its wider social and historical circumstances. It tells us very little and a lot at the same time. 

On the other hand, it is deeply radical since it pulls the rug out from under any argument that aims to foster a belief in the potential neutrality of mathematical artifacts. It means that mathematical artifacts have power. ``It is [not just] in laboratories that most new sources of power are generated''\citep[p.~163]{latour1983give}, but ever increasingly in mathematics departments and in the minds of those using and doing mathematics. Its foundational research doesn't just regularly give us the foundations for entirely new fields, such as computer science\cite{de2018making}, that shape much of our current era, but its language and tools have become the de-facto standard of most of what we understand of good (natural) science and (institutional) organisation today. Mathematical artifacts are so deeply ingrained into our everyday lives that one cannot escape them anymore, at least not if one wants to participate in society. Viewed from this perspective, they possess the same mechanisms to induce moral change that have already been explored for other technologies: they affect the making of moral decisions, our potentiality for relationships with the world and our modes of perception\cite{danaher2023mechanisms}. 
All of this now begs the question: When mathematical artifacts have politics, what are the implications for the teaching and practice of mathematics?

\section{Embedded Ethics in Mathematics}
\label{sec:embedded_ethics}
We understand the aims of a good mathematical education in accordance with classical humanistic values, i.e. the instruction of students in mathematical skills, strengthening their abilities in rational and logical reasoning, widening their curiosity, providing a sense of responsibility towards their own and other societies, people and nature, teaching them the politics, and hence ethics, of artifacts and the necessary practices to wield the power of modern mathematics for good. In essence, we want to teach them to \textit{be good people doing good mathematics} in the truest sense of the word.

This means that students need to learn to navigate the world of politics and ethics within mathematics in a gentle but intellectually stimulating way in order to foster these educational aims, to later avoid the common problem of moral overloading\cite{vandenHoven.2012}, and to be able to function at various different levels of ethical engagement\cite{chiodo2018ethics}. Such an approach includes expanding the common definition of what \textit{good} or \textit{bad} mathematics actually means\citep[p.~3]{müller2022situating}, showcasing the variety of what it can mean to be a mathematician\cite{buckmire2023definitions}, and deeply embedding ethical training into mathematics by incorporating it into the teaching at all stages of the curriculum, including appropriate mathematical exercises with a realistic and relevant ethical dimension (by going beyond ``prove this''- or ``calculate that''-type of questions and including the social, political or ethical context into them, e.g. \cite{Salado.}), in order to normalise the experience of encountering such problems\cite{chiodobursillhall}. 

It is important to remember that ``mathematicians [can be] more or less aware of philosophical tensions. But, not being philosophers (or at least not very often), they do not need to resolve these tensions — it’s enough for them to manage them by mediations and analogies''\citep[p.~19]{wagner2023paul}, and that ``ethics eludes computation because ethics is social. The concepts at the heart of ethics are not fixed or determinate in their precise meaning. To be applied they must be interpreted, and interpretations vary among individuals and groups, from context to context,
and may change over time''\citep[p.~33]{johnson2023ethical}. For approaches to embedded ethics, this implies teaching enough to navigate the issues of ethics and politics in current and future mathematical work, i.e. teaching some familiarity with philosophy is important, but it probably should not overwrite the focus on normalising discussions surrounding ethics, politics and usage, and providing the theoretical and practical knowledge necessary to navigate their often murky waters. Such foundational approaches may be particularly relevant in light of the large potential for denying social responsibility within mathematics (e.g. \cite{ernest2021mathematics, ernest2020ethics, chiodo2018ethics}) and adjacent fields (e.g. \cite{widder2023dislocated}).

\begin{figure}[ht!]
    \centering
    \includegraphics[width=0.8\textwidth]{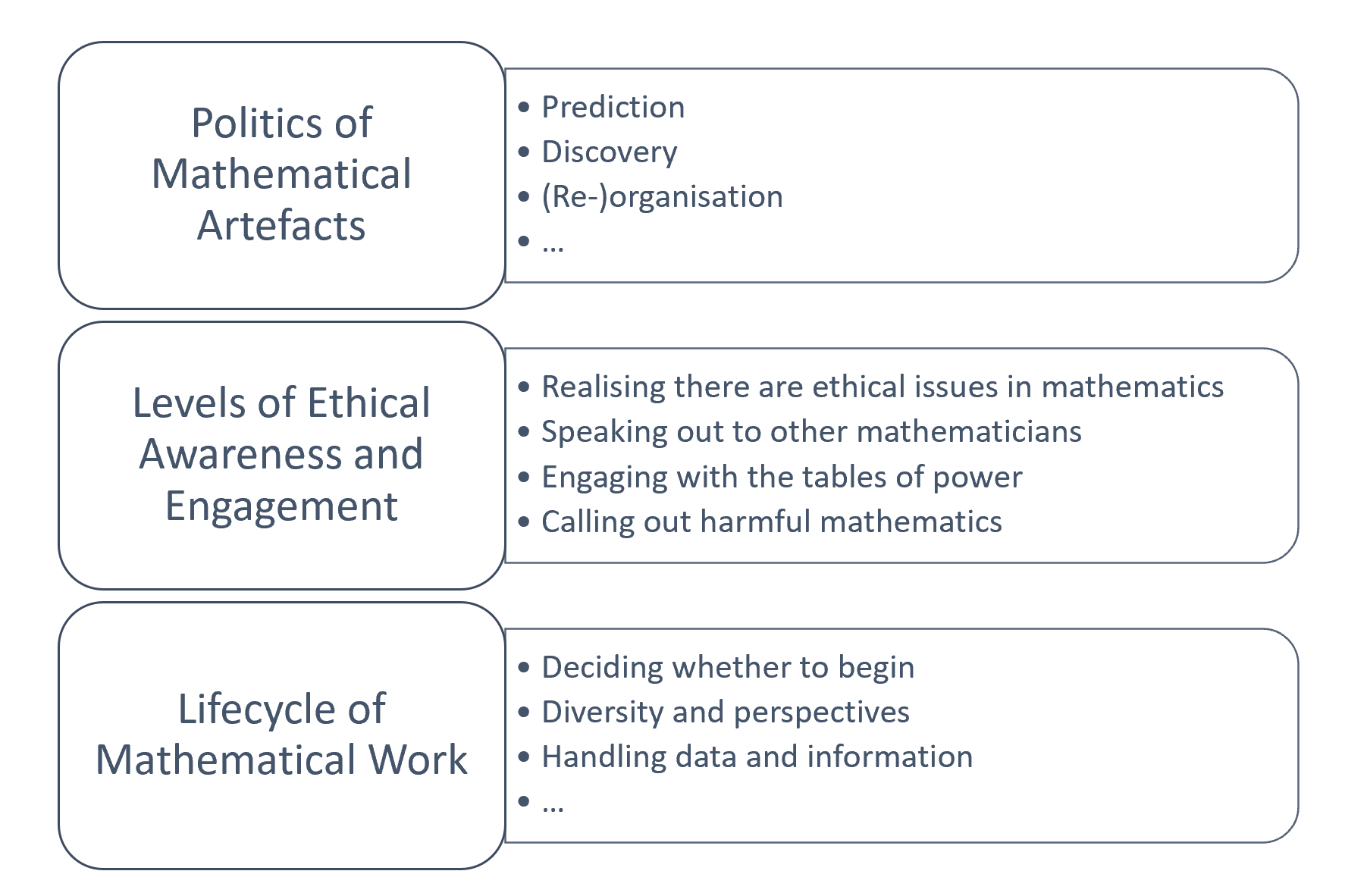}
    \caption{Figure 2: Elements of Embedded Ethics in Mathematics (based on section \ref{sec:useandeffects} and \cite{chiodo2018ethics, chiodo2023manifesto})}
    \label{fig:embedded_ethics}
\end{figure}

These insights serve as the foundation for embedding ethics into mathematical curricula. In sections \ref{sec:useandeffects} and \ref{sec:all_artifacts_have_politics} we saw the embedded politics of mathematics. But to properly transfer this into mathematical education, it needs to be combined with a knowledge of the different levels of ethical awareness that professional mathematicians and students may have\cite{chiodo2018ethics}, as well as a description or definition of the process of doing mathematics (one such description can be found in the pillars of the \textit{Manifesto for the Responsible Development of Mathematical Works}\cite{chiodo2023manifesto}). The combination of all three provides a usable way to design new mathematical exercises. We summarise this in Figure \ref{fig:embedded_ethics}.

At this point, we must note that a more hands-off approach, such as merely introducing a code of conduct or, indeed, a ``Hippocratic oath for mathematicians'' would be insufficient to accomplish these educational goals as there is no system, infrastructure and training to support (or enforce) any student or professional mathematician in their adherence to such a code \cite{muller2022hippocratic, Rittberg.2023}. It is important that any attempt to embed ethics into the mathematical curricula at universities must not just attempt to cover many of the political dimensions, but also should cover the entire process of doing mathematics in a social context (e.g. the 10 pillars outlined in the \textit{Manifesto for the Responsible Development of Mathematical Works}\cite{chiodo2023manifesto} and the different potential levels of ethical engagement\cite{chiodo2018ethics}). Only then can students see the politics of their mathematics in action.

While this task might appear daunting at first, it essentially reduces to finding the ethical or political dimensions of one's own areas of mathematical expertise if everyone is on board. Similar approaches are actively tested with good results by many universities in the subjects of computer science and artificial intelligence\cite{grosz2019embedded, bonnemains2018embedded, McLennan.2020, McLennan.2022, Cochran.2023}. It is important to remember that no single lecturer needs to design questions that cover everything, but throughout their studies, students should see as many different aspects as possible. This will automatically be assured if each lecturer covers their own courses with the help and input of others (e.g. philosophers, ethicists, practitioners, etc). The Embedded EthiCS programme at Harvard has set up a special seminar for lecturers to learn and discuss the setting of problems. Their website also includes a link to questions for many of the standard CS modules\footnote{Embedded EthiCS: Module Repository \href{https://embeddedethics.seas.harvard.edu/module}{https://embeddedethics.seas.harvard.edu/module}.}. The Cambridge University Ethics in Mathematics Project also provides useful resources for those who consider implementing embedded ethics into their programme. Their resources include exercises and additional material specifically designed for mathematicians\footnote{Cambridge University Ethics in Mathematics Project: \href{https://www.ethics.maths.cam.ac.uk/}{https://www.ethics.maths.cam.ac.uk/}.}.

\subsection{From Theory to Practice: Example Questions}
We now briefly show what exercises exploring the politics of mathematical artifacts could look like. These will be exemplary in that they show that one can construct interesting and challenging mathematical exercises that balance the mathematical difficulty and political aspects without neglecting either side. The fact that mathematical artifacts have politics does not harm mathematics, and indeed, acknowledging it helps us to create and do better mathematics. To construct these exercises, we used the insights from Figure \ref{fig:embedded_ethics} to cover at least one political dimension, one or more levels of ethical awareness, and specific pillars from the lifecycle of mathematical work.

\begin{mdframed}
\noindent\textbf{Differential Equations (Applied)}:
A detective arrives at the scene of a crime at 5:00pm. They find a warm cup of tea and measure its temperature at 40°C. By 5:30pm the tea's temperature has reduced to 30°C. 
\begin{enumerate}
    \item The police approach you with this data and ask you when the tea was likely made. Briefly discuss any questions that you still need to ask the police officers and their potential ethical relevance. What are potential barriers of communication?
\item The police are unable to provide you with more information but they ask you to give an estimate based on idealised conditions and a constant room temperature of 20°C. Giving all mathematical details and assumptions, use Newton's law of cooling to estimate when the tea was likely made.
\end{enumerate}
\end{mdframed}

This question finds its origins in the first-year differential equations course of the Mathematical Tripos at the University of Cambridge. We adjusted it to include some of the elements of embedded ethics, including
    \begin{itemize}
    \setlength\itemsep{\valitemsep}
        \item communication and deciding whether to begin (i.e. understanding mathematical assumptions and limitations of the available information), 
        \item level 1 (realising there is ethics in mathematics),
        \item as well as the political dimension of judgement, prediction and punishment.
    \end{itemize}

\begin{mdframed}\noindent\textbf{Differential Equations (Foundations)}: This question will ask you to explore the concept of dual-use within the study of differential equations.

``Dual Use Research is defined as research conducted for legitimate purposes that generates knowledge, information, technologies, and/or products that could be utilised for both benevolent and harmful purposes.'' {\small See: \href{http://www.bu.edu/research/ethics-
compliance/safety/biological-safety/ibc/dual-use-research-of-concern}{http://www.bu.edu/research/ethics-
compliance/safety/biological-safety/ibc/dual-use-research-of-concern}.}
   
\begin{enumerate}
    \item How does dual use come up in the study of differential equations?
    \item Do you know a differential equation that can be applied in a benevolent and harmful way? (Hint: Consider differential equations and their applications from your lectures. Can you apply some of them somewhere else?)
\end{enumerate}
\end{mdframed}
This question teaches students to explore the use cases of their theoretical mathematics. Students should learn the power of models and abstractions, and to experience that from this power can also come harmful use. It gently explores 
\begin{itemize}
    \setlength\itemsep{\valitemsep}
    \item level 1 (realising there is ethics in mathematics),
    \item level 4 (calling out harmful mathematics),
    \item the politics of discovery, harm, analysis and synthesis,
    \item and leads them to think about the ethics of communication, problem formulation and abstraction.
\end{itemize}

\begin{mdframed}
    \noindent\textbf{Analysis (Foundations):} This question asks you to consider what is natural about the natural numbers $\mathbb{N}$. Briefly recall the Peano axioms which we have used in our lectures to construct $\mathbb{N}$. Let $\mathbb{N}$ be a set satisfying 
    \begin{enumerate}
        \item $\mathbb{N}$ contains a special element which we call $1$.
        \item There exists a bijective map $\sigma: \mathbb{N}\rightarrow \mathbb{N} \setminus \{1\}$. 
        \item For every subset $S\subset \mathbb{N}$ such that $1\in S$ and if $n\in S$, then $\sigma(n)\in S$, it follows that $S = \mathbb{N}$.
    \end{enumerate}
    Which (if any) of $(a), (b)$, and $(c)$ feel natural to you, and why? Now consider that there are indigenous tribes who count differently. Some only have the conceptual language for the first few numbers and then use ``many'' for every larger set, e.g., they'd count $1,2,3,4, \text{ many}$. Which of $(a), (b)$ and $(c)$ would feel natural to them? Do you see why we call it ``Peano's axioms'', and not ``Peano's Theorem''? \par You can find a brief discussion of this phenomenon here: Butterworth, B. (Oct 21, 2004). What happens when you can't count past four?. \textit{The Guardian}. \url{https://www.theguardian.com/education/2004/oct/21/research.highereducation1}.
\end{mdframed}
    This question is designed to make students aware that what seems natural to them might not be natural to someone else. It uses recent research from ethnomathematics to show students that even the basics of what we perceive as pure mathematics and counting were constructed in a social context. This exercise also teaches students that Peano's axioms don't work for a finite set of numbers. It gently explores
\begin{itemize}
    \setlength\itemsep{\valitemsep}
    \item level 1 (realising there is ethics in mathematics),
    \item the politics of discovery, history and the self-image of mathematics,
    \item the role and impact of foundational axioms on how we see the world,
    \item and respect for and communication with other cultures.
\end{itemize}

\section{Conclusion}
By providing a multi-step approach, starting with specific examples showcasing that mathematical artifacts can have politics and abstracting these examples into an argument that all mathematical artifacts have politics, we constructed an applicable theory for embedded ethics in mathematics. By arguing that all mathematical artifacts have politics, we deduced the necessity to embed ethical training into mathematical curricula and showed how these insights naturally lead to new mathematical exercises. We ended the paper by briefly outlining the coverage needed for an embedded ethics curriculum in mathematics, and why this seemingly impossible task is not impossible after all. In a sense, follow Imre Lakatos, and let us all bring some well-reasoned heuristics, in this case of an ethical and political nature, into the practice of mathematics. Indeed, \textit{let us try}\citep[p.~144]{lakatos2015proofs}, even when it may appear impossible.  

\bibliographystyle{plain}
\bibliography{refs}

\end{document}